\documentclass[12pt]{amsart}
\usepackage{amsmath,amstext,amsfonts,amsthm,amssymb}
\usepackage{eucal}
\usepackage{graphicx}
\renewcommand{\subsubsection}[1]{\addtocounter{subsubsection}{1}
{\ \\[3pt]\bf \thesubsubsection. \  #1} }

\swapnumbers

{  \theoremstyle{definition}

}
\newcommand{\CT}{\operatorname{CT}}

\newcommand{\lra}{\longrightarrow}
\newcommand{\ra}{\rightarrow}

\newcommand{\tGamma}{\tilde\Gamma}

\newcommand{\tI}{\tilde I}

\newcommand{\fg}{\frak g}

\newcommand{\CI}{\mathcal{I}}

\newcommand{\BC}{\mathbb{C}}

\newcommand{\BN}{\mathbb{N}}

\newcommand{\BQ}{\mathbb{Q}}
\newcommand{\BR}{\mathbb{R}}
\newcommand{\BZ}{\mathbb{Z}}

\begin{document}


\bigskip\bigskip

\centerline{\bf REMARKS ON A TRIPLE INTEGRAL}

\bigskip\bigskip

\centerline{Bui Van Binh, Vadim Schechtman}

\vspace{1cm}

\hspace{5cm} {\it To Yuri Ivanovich Manin on his 75-th birthday}

\vspace{1cm}

\bigskip\bigskip

\centerline{\bf \S 1. Introduction}

\bigskip\bigskip

{\bf 1.1.} The following remarkable integral
$$
I_\BC(\sigma_1,\sigma_2,\sigma_3) = 
\int_{\BC^3} (1 + |x_1|^2)^{-2\sigma_1}(1 + |x_2|^2)^{-2\sigma_2}(1 + |x_3|^2)^{-2\sigma_3}
$$
$$
|x_1 - x_2|^{-2-2\nu_3}|x_2 - x_3|^{-2-2\nu_1}|x_3 - x_1|^{-2-2\nu_2} dx_1dx_2dx_3
\eqno{(1.1.1)}
$$
has appeared in [ZZ] in connection with the Liouville model of the conformal field theory. Here $\sigma_i\in \BC$, 
$$
\nu_1 = \sigma_1 - \sigma_2 - \sigma_3,\  
\nu_2 = \sigma_2 - \sigma_3 - \sigma_1,\ 
\nu_3 = \sigma_3 - \sigma_1 - \sigma_2
$$
et $dx$ denotes the standard Haar measure on $\BC$. 

Set
$$
\tI_\BC(\sigma_1,\sigma_2,\sigma_3) = \int_{\BC^2} (1 + |x_1|^2)^{-2\sigma_1}(1 + |x_2|^2)^{-2\sigma_2}
|x_1 - x_2|^{-2-2\nu_3} dx_1dx_2
\eqno{(1.1.2)}
$$
To compute (1.1.1) the authors of [ZZ] first note that 
$$
I_\BC(\sigma_1,\sigma_2,\sigma_3) = \pi \tI_\BC(\sigma_1,\sigma_2,\sigma_3)
\eqno{(1.1.3)}
$$
This may be proven using an $SU(2)$-symmetry (cf. 5.1 below). So (1.1.1) and (1.1.2) converge for 
$\Re\sigma_1, \Re\sigma_2, \Re\nu_3$ sufficiently large.  

Then the authors give (without proof) the value of $\tI_\BC$ and hence that of $I_\BC$: 
$$
I_\BC(\sigma_1,\sigma_2,\sigma_3) = \pi^3\frac{\Gamma(\sigma_1 + \sigma_2 + \sigma_3 - 1)
\Gamma(-\nu_1)\Gamma(-\nu_2)\Gamma(-\nu_3)}{\Gamma(2\sigma_1)\Gamma(2\sigma_2)\Gamma(2\sigma_3)}
\eqno{(1.1.4)}
$$
A proof (somewhat artificial) of (1.1.4) may be found in [HMW]. We propose another proof in 5.3 below. 



{\bf 1.2.} In this note we take the study of the {\it real}, {\it $q$-deformed} and {\it $p$-adic} versions of (1.1.1). 

A real version of (1.1.1) is the integral (1.2.1) below, cf. \S 6.  
It has appeared in [BR] in connection with a study of periods of automorphic triple products: 
$$
I_\BR(\sigma_1,\sigma_2,\sigma_3) := \frac{1}{(2\pi)^3}\int_{-\pi}^\pi \int_{-\pi}^\pi \int_{-\pi}^\pi\ \prod_{i=1}^3 |\sin(\theta_i - \theta_{i+1})|^{(\nu_{i+2}-1)/2} d\theta_1 d\theta_2 d\theta_3 =  
$$
$$
= \frac{\Gamma((\nu_1 + 1)/4)\Gamma((\nu_2 + 1)/4)\Gamma((\nu_3 + 1)/4)\Gamma((\sum_i\nu_i + 1)/4)}
{\Gamma(1/2)^3\Gamma((1 - \sigma_1)/2)\Gamma((1 - \sigma_2)/2)\Gamma((1 - \sigma_3)/2)}
\eqno{(1.2.1)}
$$
The index $i$ under the integral is understood modulo $3$. The authors of [BR] provide an elegant 
proof of (1.2.1) using Gaussian integrals. 

In this note we propose and calculate a $q$-deformation 
of this integral (see Thm 1.6 and \S 4 below); in the limit $q\ra 1$ this gives (1.2.1). 

Set 
$$
a_i = \frac{\nu_i - 1}{4}
$$
Then (1.2.1) rewrites as 
$$
J(a_1,a_2,a_3) := \frac{1}{(2\pi)^3}\int_{-\pi}^{\pi} \int_{-\pi}^{\pi} \int_{-\pi}^{\pi}\ \prod_{i=1}^3 |\sin(\theta_i - \theta_{i+1})|^{2a_{i+2}} d\theta_1 d\theta_2 d\theta_3 =  
$$
$$
= \frac{\Gamma(a_1+1/2)\Gamma(a_2+1/2)\Gamma(a_3+1/2)\Gamma(\sum a_i + 1)}
{\Gamma(1/2)^3\Gamma(a_1 + a_2 + 1)\Gamma(a_2 + a_3 + 1)\Gamma(a_3 + a_1 + 1)}
$$
$$
= \frac{\prod_i \Gamma(2a_i + 1)\Gamma(\sum a_i + 1)}
{4^{\sum a_i}\prod_i \Gamma(a_i + 1)\prod_{i<j} \Gamma(a_i + a_j + 1)}
\eqno{(1.2.2)}
$$
where we have used the duplication formula
$$
\Gamma(2a+1) = 2^{2a}\pi^{-1/2}\Gamma(a+1/2)\Gamma(a+1)
$$
Let us suppose that $a_i$ are positive integers. After a change of variables $y_j = e^{2i\theta_j}$ it is easy 
to see that (1.2.2) is equivalent to 
$$
\CT\prod_{1\leq i < j\leq 3}(1 - y_i/y_j)^{a_{ij}}(1 - y_j/y_i)^{a_{ij}}  
$$
$$
= \frac{(a_1+a_2+a_3)!\prod_{i=1}^3 (2a_i)!}{\prod_{i=1}^3 a_i!\prod_{1\leq i< j\leq 3}(a_i+a_j)!}
\eqno{(1.2.3)}
$$
Here $a_{ij} := a_k$ where $\{k\} = \{1,2,3\}\setminus\{i,j\}$ and $\CT$ means "constant term".  

If $a_1 = a_2 = a_3 = a$, this gets into   
$$
\CT\prod_{1\leq i \neq j \leq 3} (1 - y_i/y_j)^{a} = 
\frac{\Gamma(3a + 1)}
{\Gamma(a + 1)^3}
\eqno{(1.2.4)}
$$
which is the classical Dyson formula for the root system $A_2$, cf. [Dy], (142). 

In \S 3 below we give an independent proof of (1.2.3). 

{\bf 1.3.} The reformulation (1.2.3) allows us to write down a $q$-deformation of it. It looks 
as follows. For a natural $a$ denote as usual 
$$
(x;q)_a = \prod_{i=0}^{a-1} (1 - xq^i) = \frac{(x;q)_\infty}{(xq^a;q)_\infty}
$$
where
$$
(x;q)_\infty =  \prod_{i=0}^{\infty} (1 - xq^i)
\eqno{(1.3.1)}
$$
Here $q$ is a formal variable. 

Denote
$$
[a]^!_q = \frac{(q;q)_a}{(1 - q)^a}
$$

{\bf 1.4. Theorem.} {\it Let $a_1, a_2, a_3$ be natural. Then
$$
\CT\prod_{1\leq i < j\leq 3} (y_i/y_j;q)_{a_{ij}}(qy_j/y_i;q)_{a_{ij}}   
$$
$$
= \frac{[a_1+a_2+a_3]_q^!\prod_{i=1}^3 [2a_i]_q^!}{\prod_{i=1}^3 [a_i]_q^!\prod_{1\leq i< j\leq 3}[a_i+a_j]_q^!}
\eqno{(1.4.1)}
$$}

For a proof see \S 4 below. 

If $a_1 = a_2 = a_3$  this becomes the (proven) Macdonald's $q$-constant term conjecture for 
the root system $A_2$. 

In fact, (1.4.1) is in turn a particular case of the following beautiful formula due to W.Morris. 

{\bf Theorem} (W.Morris). {\it Let $a_1, a_2, a_3$ be natural and $\sigma\in \Sigma_3$ be an arbitrary 
permutation. Then 
$$
\CT\prod_{1\leq i < j\leq 3} (y_i/y_j;q)_{a_{ij}}(qy_j/y_i;q)_{a_{\sigma(i)\sigma(j)}}   
$$
$$
= \frac{[a_1+a_2+a_3]_q^!\prod_{i=1}^3 [a_i + a_{\sigma(i)}]_q^!}{\prod_{i=1}^3 [a_i]_q^!\prod_{1\leq i< j\leq 3}[a_i+a_j]_q^!}
\eqno{(1.4.2)}
$$}

This is the case of the {\it $A_2$-isolated labeling} of Morris' conjecture, cf. [Mo], 4.3. 
The proof of (1.4.2) is contained in {\it op. cit.} 5.12. 
The formula (1.4.1) is the case $\sigma = $ the identity permutation of (1.4.2). 
Our proof of (1.4.1) is less involved than Morris' proof of the general case.

\bigskip

{\bf 1.5.} Let us generalize Thm. 1.4 to the case of complex $a_i$. To this end we suppose 
that $q$ is a real number, $0 < q < 1$. 

For any $x, a\in \BC$ we define 
$$
(x;q)_a = \frac{(x;q)_\infty}{(xq^a;q)_\infty}
$$
Define as usually
$$
\Gamma_q(x) = (1 - q)^{1-x}\frac{(q,q)_\infty}{(q^x,q)_\infty}
$$

We denote by $T^3$ the torus
$$
T^3 = \{(y_1,y_2,y_3)\in \BC^3|\ |y_i| = 1,\ 1\leq i\leq 3\}
$$

{\bf 1.6. Theorem.} {\it For $\Re(a_i) > 0,\ 1\leq i \leq 3$ 
$$
\frac{1}{(2\pi i)^3}\int_{T^3} \prod_{1\leq i < j\leq 3} (y_i/y_j;q)_{a_{ij}}(qy_j/y_i;q)_{a_{ij}} 
\frac{dy_1}{y_1}\frac{dy_2}{y_2}\frac{dy_3}{y_3} = 
$$
$$
\frac{\Gamma_q(a_1+a_2+a_3+1)\prod_{i=1}^3 \Gamma_q(2a_i+1)}{\prod_{i=1}^3 \Gamma_q(a_i+1)\prod_{1\leq i< j\leq 3}\Gamma_q(a_i+a_j+1)}
\eqno{(1.6.1)}
$$}

For a proof see 4.6 below. Passing to the limit $q\ra 1$ gives (1.2.1).

{\bf 1.7.} Let us describe a $p$-adic version of (1.1.1). 

Let $p$ be a prime number; consider the field $\BQ_p$ of rational $p$-adic numbers. 
Let $d_px$ denote the Haar measure on $\BQ_p$ normalized by the 
condition 
$$
\int_{\BZ_p} d_px = 1
$$
Let
$$
|.|_p:\ \BQ_p^\times \lra \BR^\times_{> 0}
$$
be the standard $p$ adic norm, $|p|_p = p^{-1}$; we set $|0|_p = 0$. 

We have $d_p(ax) = |a|_p d_px$, so $|a|_p$ is a $p$-adic analog of $|z|^2, z\in \BC$ (sic!). 

Define a function $\psi_p(x), x\in \BQ_p$, by 
$$
\psi_p(x) = \max\{|x|_p, 1\}
\eqno{(1.7.1)}
$$
This is an analog of $|z|^2 + 1, z\in \BC$, see 2.5 below. 

Set
$$
\Gamma_{\BQ_p}(\sigma) = \frac{1 - p^{-1}}{1 - p^{-\sigma}},\ \sigma\in \BC
$$

{\bf 1.8.} The following integrals are $p$-adic analogs of (1.1.1):
$$
I_{\BQ_p}(\sigma_1,\sigma_2,\sigma_3) = \int_{\BQ_p^3} \psi_p(x_1)^{-2\sigma_1}\psi_p(x_2)^{-2\sigma_2}\psi_p(x_3)^{-2\sigma_3}
$$
$$
|x_1 - x_2|_p^{-1-\nu_3}|x_2 - x_3|_p^{-1-\nu_1}|x_3 - x_1|_p^{-1-\nu_2} d_px_1d_px_2d_px_3
\eqno{(1.8.1)}
$$
and of (1.1.2):
$$
\tI_{\BQ_p}(\sigma_1,\sigma_2,\sigma_3) = \int_{\BC^2} \psi_p(x_1)^{-2\sigma_1}\psi_p(x_2)^{-2\sigma_2}
|x_1 - x_2|_p^{-1-\nu_3} d_px_1d_px_2
\eqno{(1.8.2)}
$$ 

{\bf 1.9. Theorem.} (i)
$$
I_{\BQ_p}(\sigma_1,\sigma_2,\sigma_3) = \frac{1}{\Gamma_p(2)}\tI_{\BQ_p}(\sigma_1,\sigma_2,\sigma_3)
\eqno{(1.9.1)}
$$
(ii)
$$
\tI_{\BQ_p}(\sigma_1,\sigma_2,\sigma_3) =  \frac{\Gamma_{\BQ_p}(\sigma_1 + \sigma_2 + \sigma_3 - 1)
\Gamma_{\BQ_p}(-\nu_1)\Gamma_{\BQ_p}(-\nu_2)\Gamma_{\BQ_p}(-\nu_3)}{\Gamma_{\BQ_p}(2\sigma_1)\Gamma_{\BQ_p}(2\sigma_2)\Gamma_{\BQ_p}(2\sigma_3)}
\eqno{(1.9.2)}
$$

For a proof, see \S 2 below. 

{\bf 1.10.} Let $G = PGL(2)$. The integral $I_\BC$ (resp. $I_{\BR}, I_{\BQ_p}$) is related to  invariant functionals 
on triple products $V_1\otimes V_2\otimes V_3$ where $V_i$ are irreducible $G(K)$-representations of the principal 
series, with $K = \BC$ (resp. $\BR$ or $\BQ_p$), cf. [BR] for the real case. So its $q$-deformation (1.6.1) should be related to 
the same objects connected with the quantum group $U_q\fg(\BR)$ where $\fg = Lie\ G$.

{\bf 1.11. Notation.} $\BN = \{0, 1, 2, \ldots\}$.      

{\bf 1.12. Acknowledgement.} The second author is thankful to J.Bernstein and A.Reznikov for inspiring discussions. 
We thank W.Zudilin for sending the unpublished dissertation of W.Morris. 

\bigskip

\bigskip\bigskip

\centerline{\bf \S 2. The $p$-adic case}

\bigskip\bigskip

{\bf 2.1. Notation.} Let
$$
v_p:\ \BQ_p \lra \BZ\cup\{\infty\}
$$
denote the usual $p$-adic valuation, i.e. $v_p(x) = n$ if $x\in p^n\BZ_p\setminus p^{n+1}\BZ_p$, 
$v_p(0) = \infty$. 

For $n, m\in \BZ$ we denote
$$
A_{\leq n} = \{x\in \BQ_p|\ v_p(x) \leq n\},\ 
A_{\geq n} = \{x\in \BQ_p|\ v_p(x) \geq n\},\ 
$$
$$ 
A_{[n,m]} = A_{\geq n}\cap A_{\leq m},\ A_n = A_{[n,n]}
$$

We also set
$$
\Gamma_{\BQ_p}(\infty) := \lim_{a\ra\infty} \Gamma_{\BQ_p}(a) = 1 - p^{-1}
$$

{\bf 2.2. A $p$-adic hypergeometric function.} Define
$$
F_{\BQ_p}(a,c;y) = \int_{\BQ_p} \psi_p(x)^a |x - y|_p^c d_px,
$$
$a, c \in \BC; y\in \BQ_p$.

{\bf 2.3. Lemma.} (i) {\it If $v_p(y) \geq 0$ then
$$
F_{\BQ_p}(a,c;y) = \Gamma_{\BQ_p}(c+1) - \Gamma_{\BQ_p}(a+c+1)
$$}

(ii) {\it If $v_p(y) = n < 0$ then
$$
F_{\BQ_p}(a,c;y) = 
p^{-n(a+c+1)}\Gamma_{\BQ_p}(c+1) - p^{-n(a+c+1)}\Gamma_{\BQ_p}(a+c+1)
$$
$$
+ \frac{p^{-nc}\Gamma_{\BQ_p}(\infty)}{\Gamma_{\BQ_p}(n(a+1)+1)}
- \frac{p^{-nc}\Gamma_{\BQ_p}(\infty)\Gamma_{\BQ_p}(a+1)}{\Gamma_{\BQ_p}(n(a+1))}
$$}

{\bf Proof.} Let us denote for brevity
$$
f(a,c;x,y) = \psi_p(x)^a |x - y|_p^c
$$ 

(i) Let $v_p(y) = n \geq 0$. Decompose $\BQ_p$ into the following 
areas:
$$
\BQ_p = A_{< 0}\cup A_{[0,n-1]}\cup A_n\cup A_{\geq n} 
$$
Then
$$
\int_{A_{< 0}} f d_px = - \Gamma_{\BQ_p}(a+c+1),
$$
$$
\int_{A_{[0,n-1]}} f d_px = (1 - p^{-n(c+1)})\Gamma_{\BQ_p}(c+1),
$$
$$
\int_{A_{>n}} f d_px = p^{-n(c+1)-1}
$$
To evaluate $\int_{A_n} f(x,y) d_px$, we decompose $A_n$ into two subsets depending on $y\in A_n$: 
$A_n = A'_n(y) \cup A_n''(y)$ where
$$
A'_n(y) = \{x\in A_n|\ v_p(x - y) = n\},\ A''_n(y) = \{x\in A_n|\ v_p(x - y) > n\}
\eqno{(2.3.1)}
$$
Then
$$
\int_{A'_n(y)} f d_px = (p-2)p^{-n(c+1)-1}
$$
and
$$
\int_{A''_n(y)} f d_px = p^{-(n+1)(c+1)}\Gamma_{\BQ_p}(c+1),
$$
so that
$$
\int_{A_n} f d_px = (p-2)p^{-n(c+1)-1} + p^{-(n+1)(c+1)}\Gamma_{\BQ_p}(c+1). 
$$
Adding up, we get (i). 

(ii) is proved in a similar manner. Let $v_p(y) = n < 0$.  We decompose
$$
\BQ_p = A_{< n}\cup A_n\cup A_{[-n+1,-1]}\cup A_{\geq 0},
$$
and for $y\in A_n$
$$
A_n = A'_n(y)\cup A''_n(y)
$$
as in (2.3.1). Then
$$
\int_{A_{< n}} f d_px = - p^{-n(a+c+1)}\Gamma_{\BQ_p}(a+c+1),
$$
$$
\int_{A_{[-n+1,-1]}} f d_px = (1 - p^{-1})p^{-nc}\sum_{m=n+1}^{-1} p^{-m(a+1)},
$$
$$
\int_{A_{\geq 0}} f d_px = p^{-nc},
$$
$$
\int_{A'_n(y)} f d_px = (p - 2)p^{-n(a+c+1)-1}
$$
and
$$
\int_{A''_n(y)} f d_px = p^{-(n+1)(c+1)}\Gamma_{\BQ_p}(c+1)
$$
Adding up, we get (ii). $\square$. 

{\bf 2.4. Theorem.} 
$$
J(a,b,c) := \int\int_{\BQ_p^2} \psi_p(x)^a\psi_p(y)^b |x - y|_p^c d_px d_py = 
$$
$$
= \frac{\Gamma_{\BQ_p}(c+1)\Gamma_{\BQ_p}(-a-c-1)\Gamma_{\BQ_p}(-b-c-1)\Gamma_{\BQ_p}(-a-b-c-2)}
{\Gamma_{\BQ_p}(-a)\Gamma_{\BQ_p}(-b)\Gamma_{\BQ_p}(-a-b-2c-2)}
\eqno{(2.4.1)}
$$

{\bf Proof.} By definition
$$
J(a, b, c) = \int_{\BQ_p} \psi_p(y)^b F_{\BQ_p}(a,c;y) d_py
$$
Using Lemma 2.3 we readily compute this integral and arrive at (2.4.1). $\square$

This theorem is equivalent to  (1.9.2).

{\bf 2.5. Proof of} (1.9.1).  We shall use the same method as in the complex case, cf. \S 5 below.

Let
$$
K = SL_2(\BZ_p) \subset G = SL_2(\BQ_p)
$$
If $v\in \BQ_p^2$ and $g\in K$ then 
$$
|v|_p = |gv|_p
\eqno{(2.5.1)}
$$
where 
$$
|(a, b)|_p = \max\{|a|_p,|b|_p\}
$$
For 
$$
g = \left(\begin{matrix} a & b\\ c & d\end{matrix}\right)\in G,\ z\in \BQ_p
$$ 
set
$$
g\cdot z = \frac{az + b}{cz + d}
$$
It follows:
$$
\psi_p(g\cdot x) = \frac{\psi_p(x)}{|cx + d|_p},\ g\in K
\eqno{(2.5.2)}
$$
We have also for $g\in G$
$$
g\cdot x - g\cdot y = \frac{x - y}{(cx + d)(cy + d)}
\eqno{(2.5.3)}
$$
and
$$
d_p(g\cdot z) = \frac{d_pz}{|cz + d|_p^2},
\eqno{(2.5.4)}
$$
cf. [GGPS], Ch. II, \S 3, no. 1. 

We have 
$$
I_{\BQ_p}(\sigma_1,\sigma_2,\sigma_3) = \int_{\BQ_p} \psi_p(x_3)^{-2\sigma_3}\biggl(\int_{\BQ_p^2}\psi_p(x_1)^{-2\sigma_1}\psi_p(x_2)^{-2\sigma_2}
$$
$$
|x_1 - x_2|_p^{-1-\nu_3}|x_2 - x_3|^{-1-\nu_1}|x_3 - x_1|^{-1-\nu_2} d_px_1d_px_2\biggr)d_px_3
\eqno{(2.5.5)}
$$
Given $y\in \BQ_p$, set $a(y) = y$ if $y\in \BZ_p$ and $a(y) = y^{-1}$ if $y\notin\BZ_p$; 
so $a(y)\in \BZ_p$ in any case. 

Define a matrix $k(y)\in K$ by:

(i) if $y\in \BZ_p$ then
$$
k(y) = \left(\begin{matrix} 1 & -a(y)\\ 0 & 1\end{matrix}\right)
$$

(ii) if $y\in\BQ_p\setminus \BZ_p$ then
$$
k(y) = \left(\begin{matrix} a(y) & -1\\ 1 & 0\end{matrix}\right)
$$
In the internal integral in (2.5.5) let us make a change of variables 
$$
x_i = k(x_3)^{-1}\cdot y_i,\ 
i = 1, 2
$$
Using (2.5.2) - (2.5.4) we get
$$
I_{\BQ_p}(\sigma_1,\sigma_2,\sigma_3) = \int_{\BQ_p} \frac{d_px_3}{\psi_p(x_3)^2}\cdot I'_{\BQ_p}(\sigma_1,\sigma_2,\sigma_3)
$$ 
where
$$
\cdot I'_{\BQ_p}(\sigma_1,\sigma_2,\sigma_3) = 
$$
$$
\int\int_{\BQ_p^2} |y_1|_p^{-1-\nu_2}|y_2|_p^{-1-\nu_1}
\psi_p(y_1)^{-2\sigma_1}\psi_p(y_2)^{-2\sigma_2}|y_1 - y_2|_p^{-1-\nu_3} d_py_1d_py_2
$$
After one more substitution $y_i \mapsto y_i^{-1},\ i = 1, 2$, 
$$
I'_{\BQ_p}(\sigma_1,\sigma_2,\sigma_3) = \tI_{\BQ_p}(\sigma_1,\sigma_2,\sigma_3)
$$
(note that $d_p(y^{-1}) = d_py/|y|_p^2$). 
Finally we conclude by the following easily proved $p$-adic version of (5.1.6):

{\bf 2.6. Lemma.}
$$
\int_{\BQ_p} \frac{d_px}{\psi_p(x)^2} = 1 + p^{-1}
$$
$\square$

\bigskip\bigskip



\centerline{\bf \S 3. The real case}

\bigskip\bigskip

{\bf 3.1. Theorem.} {\it If $a, b, c\in \BN$ then 
$$
\CT (1 - y_1/y_2)^c(1 - y_2/y_1)^c(1 - y_1/y_3)^b(1 - y_3/y_1)^b(1 - y_2/y_3)^a(1 - y_3/y_2)^a
$$
$$
= \frac{(2a)!(2b)!(2c)!(a + b + c)!}{a!b!c!(a + b)!(a + c)!(b + c)!}
$$}

{\bf 3.2. Lemma} (A.C.Dixon's identity). 
$$
\sum_{n = -\infty}^\infty\ (-1)^n\binom{a+b}{a+n}\binom{b+c}{b+n}\binom{a+c}{c+n}
= \frac{(a + b + c)!}{a!b!c!}
$$

(We set $\binom{a}{b} = 0$ for $b < 0$.) 

See [K], 1.2.6, Exercice 62 and Answer to Ex. 62, p. 490 (one finds also interesting references there). $\square$

{\bf 3.3. Proof of 3.1.} The Laurent polynomial on the left hand side is
$$
f(y_1,y_2,y_3) = \sum_{i,j,k} (-1)^{a+b+c-i-j-k} 
\binom{2c}{i}\binom{2a}{j}\binom{2b}{k}y_1^{b-c+i-k}y_2^{c-a+j-i}y_3^{a-b+k-j}
$$
whence the constant term corresponds to the values
$$
b-c+i-k = c-a+j-i = a-b+k-j = 0
$$
Set $n = c - i$; then $n = b - k = a - j$ as well, so
$$
\CT f(y_1,y_2,y_3) = \sum_{n=0}^\infty (-1)^{-3n}\binom{2a}{a+n}\binom{2b}{b+n}\binom{2c}{c+n} = 
$$ 
$$
\frac{(2a)!(2b)!(2c)!}{(a + b)!(a + c)!(b + c)!}\sum_{n=0}^\infty\ (-1)^n\binom{a+b}{a+n}\binom{b+c}{b+n}\binom{a+c}{c+n},
$$
and the application of 3.2 finishes the proof. $\square$

\bigskip\bigskip



\centerline{\bf \S 4. The real $q$-deformed case}

\bigskip\bigskip

{\bf 4.1.} Notation:
$$
\left[\begin{matrix} a\\ b\end{matrix}\right]_q = \frac{[a]^!_q}{[b]^!_q[a-b]^!_q}
$$
If $a, b\in \BZ, a\geq 0, b < 0$ we set
$$
\left[\begin{matrix} a\\ b\end{matrix}\right]_q = 0
$$

Set
$$
u = x_2/x_1,\ v = x_3/x_2,\ w = x_1/qx_3;
$$
$$
a = a_3,\ b = a_1,\ c = a_2.
$$
We are interested in the constant term of
$$
F_q(u, v, w) = (qu;q)_a(u^{-1};q)_a(qv;q)_b(v^{-1};q)_b(qw;q)_b(w^{-1};q)_b
\eqno{(4.1.1)}
$$
where $uvw = q^{-1}$.

{\bf 4.2. Lemma} (K.Kadell).  
$$
(qx;q)_b(x^{-1};q)_a = \sum_{i = -a}^b\ q^{i(i+1)/2}\left[\begin{matrix} a+b\\ a+i\end{matrix}\right]_q(-x)^i
$$

See [Ka], (3.31). $\square$  

{\bf 4.3. Lemma} ($q$-Dixon identity). {\it On has two equivalent formulas:} 

(i) 
$$
\sum_{n = -\infty}^\infty\ (-1)^n q^{n(3n+1)/2}\left[\begin{matrix} a + b\\ a + n\end{matrix}\right]_q
\left[\begin{matrix} b + c\\ b + n\end{matrix}\right]_q\left[\begin{matrix} c + a\\ c + n\end{matrix}\right]_q =
$$
$$
= \frac{[a + b + c]_q^!}{[a]_q^![b]_q^![c]_q^!}
\eqno{(4.3.1)}
$$

(ii)
$$
\sum_{n = -\infty}^\infty\ (-1)^n q^{n(3n+1)/2}\left[\begin{matrix} 2a\\ a + n\end{matrix}\right]_q
\left[\begin{matrix} 2b\\ b + n\end{matrix}\right]_q\left[\begin{matrix} 2c\\ c + n\end{matrix}\right]_q =  
$$
$$
= \frac{[2a]_q^![2b]_q^![2c]_q^![a + b + c]_q^!}{[a]_q^![b]_q^![c]_q^![a+b]_q^![b+c]_q^![a+c]_q^!}
\eqno{(4.3.2)}
$$  

Cf. [K], answer to Exercice 1.2.6, [C]. $\square$

{\bf 4.4.} Now we can prove (1.4.1). Replace in the product (4.1.1) the double products 
like $(qu;q)_a(u^{-1};q)_a$ using 4.2. In the resulting expression the contant term will be the 
sum of coefficients at $u^iv^iv^i$ divided by $q^i$. Thus
$$
\CT F_q(u, v, w) = \sum_{n = -\infty}^\infty\ (-1)^n q^{-n}q^{3n(n+1)/2}\left[\begin{matrix} 2a\\ a + n\end{matrix}\right]_q
\left[\begin{matrix} 2b\\ b + n\end{matrix}\right]_q\left[\begin{matrix} 2c\\ c + n\end{matrix}\right]_q 
$$
$$
= 
\sum_{n = -\infty}^\infty\ (-1)^n q^{n(3n+1)/2}\left[\begin{matrix} 2a\\ a + n\end{matrix}\right]_q
\left[\begin{matrix} 2b\\ b + n\end{matrix}\right]_q\left[\begin{matrix} 2c\\ c + n\end{matrix}\right]_q 
$$
$$
= \frac{[2a]_q^![2b]_q^![2c]_q^![a + b + c]_q^!}{[a]_q^![b]_q^![c]_q^![a+b]_q^![b+c]_q^![a+c]_q^!}
$$
by (4.3.2). This finishes the proof of (1.4.1). $\square$

Let us prove Thm. 1.6. We shall use an idea going back to Hardy, cf. [B], 5.5; [S]; [M], 17.2. 

First (1.6.1) is true if all $a_i\in \BN$ --- this is Thm. 1.4. 
Now we shall use 

{\bf 4.5. Lemma.} {\it Let $f(z)$ be a function holomorphic and bounded for $\Re z \geq 0$ such that 
$f(z) = 0$ for $z\in \BN$. Then $f(z) \equiv 0$.}

This is a particular case of  {\bf Carlsson's theorem}, cf. [B], 5.3; [T], 5.8.1. $\square$

{\bf 4.6.} Set
$$
\tGamma_q(a) = \frac{\prod_{i=1}^\infty (1 - q^i)}{\prod_{i=0}^\infty (1 - q^{a+i})}, 
$$
so that 
$$
\Gamma_q(a) = (1 - q)^{1 - a}\tGamma_q(a);
$$
Recall that $0 < q < 1$. 

We have 
$$
1 - |b|q^{\Re s} \leq |1 - bq^s|\leq 1 + |b|q^{\Re s},\ \Re s \geq 0,
$$
and 
$$
1 + t \leq e^t,\ t\geq 0.
$$
It follows: 
$$
|\prod_{i=0}^\infty (1 - q^{a+i})|\leq \prod_{i=0}^\infty (1 + q^{\Re a + i}) \leq 
$$
$$
\prod_{i=0}^\infty e^{q^{\Re a + i}} = e^{\sum_{i=0}^\infty q^{\Re a + i}} = 
e^{q^{\Re a}/(1 - q)}\leq e^{1/(1 - q)}
$$
for $\Re a \geq 0$. 

On the other hand
$$
|\prod_{i=0}^\infty (1 - q^{a+i})|\geq \prod_{i=0}^\infty (1 - q^{\Re a + i}) 
\geq \prod_{i=0}^\infty (1 - q^{a_0 + i})
$$
for $\Re a \geq a_0 > 0$. 

Fix $a_0 > 0$. It follows that there exist constants $C_1, C_2 > 0$ such that 
$$
C_1 \leq \tGamma_q(a) \leq C_2
\eqno{(4.6.1)}
$$
for all $a,\ \Re a \geq a_0$.  

Consider the right hand side of (1.6.1)
$$
f(a_1, a_2, a_3) := 
\frac{\Gamma_q(a_1+a_2+a_3+1)\prod_{i=1}^3 \Gamma_q(2a_i+1)}
{\prod_{i=1}^3 \Gamma_q(a_i+1)\prod_{1\leq i< j\leq 3}\Gamma_q(a_i+a_j+1)} =
$$
$$
\frac{\tGamma_q(a_1+a_2+a_3+1)\prod_{i=1}^3 \tGamma_q(2a_i+1)}
{\prod_{i=1}^3 \tGamma_q(a_i+1)\prod_{1\leq i< j\leq 3}\tGamma_q(a_i+a_j+1)}
$$
It follows from (4.7.1) that there a constant $C_3 > 0$ such that 
$$
|f(a_1, a_2, a_3)| \leq C_3
$$
for all $a_1, a_2, a_3$ with the real part $\geq a_0$.  

In the same manner we prove that if $g(a_1, a_2, a_3;x_1, x_2, x_3)$ is the expression 
under the integral from the left hand side of (1.6.1), there exist a constant $C_4 > 0$ such that
$$
|g(a_1, a_2, a_3)| \leq C_4
$$
for all $a_1, a_2, a_3$ with the real part $\geq a_0$ and $(x_1, x_2, x_3)\in T^3$; thus 
$$
h(a_1, a_2, a_3) = \frac{1}{(2\pi^3)}|\int_{T^3} g(a_1, a_2, a_3;x_1, x_2, x_3) dx_1dx_2dx_3|
$$
is also bounded by a constant not depending on $a_i$.

By Thm. 1.4 we know that $h(a_1, a_2, a_3) = f(a_1, a_2, a_3)$ if all $a_i\in \BN$. 
Now applying ($3$ times) 4.5 we conclude that this is true for all $a_i$ with 
$\Re a_i \geq a_0$. This proves Thm. 1.6. $\square$

\bigskip\bigskip



\centerline{\bf \S 5. The complex case}

\bigskip\bigskip

{\bf 5.1. Proof of} (1.1.3). We give some details because we use exactly the same argument in the $p$-adic 
case, cf. 2.5. 

Let 
$$
K = SU(2) = \{\left(\begin{matrix} a & b\\ -\bar b & \bar a\end{matrix}\right)|\ a, b\in \BC,\ |a|^2 + |b|^2 = 1\}\subset
$$
$$ 
\subset G = SL_2(\BC)
\eqno{(5.1.1)}
$$

If $v\in \BC^2$ and $g\in K$ then 
$$
|v| = |gv|
\eqno{(5.1.2)}
$$
where 
$$
|(a, b)|^2 = |a|^2 + |b|^2
$$
For 
$$
g = \left(\begin{matrix} a & b\\ c & d\end{matrix}\right)\in G,\ z\in \BC
$$ 
set
$$
g\cdot z = \frac{az + b}{cz + d}
$$
It follows:
$$
1 + |g\cdot x|^2 = \frac{1 + |x|^2}{|cx + d|^2},\ g\in K
\eqno{(5.1.3)}
$$
We have also for $g\in G$
$$
g\cdot x - g\cdot y = \frac{x - y}{(cx + d)(cy + d)}
\eqno{(5.1.4)}
$$
and
$$
d(g\cdot z) = \frac{dz}{|cz + d|^4}
\eqno{(5.1.5)}
$$
(recall that in real coordinates $z = x + iy$ we have $dz = dx dy$).  

Using this, let us evaluate
$$
I_\BC(\sigma_1,\sigma_2,\sigma_3) = 
$$
$$
\int_\BC dx_3 (1 + |x_3|^2)^{-2\sigma_3}\biggl( 
\int\int_{\BC^2} \prod_{i=1}^2\ (1 + |x_i|^2)^{-2\sigma_i} 
\prod_{i=1}^3 |x_i - x_{i+1}|^{-2-2\nu_{i+2}} 
dx_1dx_2\biggr)
$$

In the internal integral let us make a change of variables $x_i = k^{-1}\cdot y_i,\ 
i = 1, 2$ with $k\in K$ as in (5.1.1) with 
$$
a = \frac{1}{(1 + |x_3|^2)^{1/2}},\ b = -  \frac{x_3}{(1 + |x_3|^2)^{1/2}},
$$
so $y_3 = k\cdot x_3 = 0$. 

We get
$$
I_\BC(\sigma_1,\sigma_2,\sigma_3) = \int_\BC \frac{dx_3}{(1 + |x_3|^2)^2}\cdot I'_\BC(\sigma_1,\sigma_2,\sigma_3)
$$ 
where
$$
I'_\BC(\sigma_1,\sigma_2,\sigma_3) = 
$$
$$
\int\int_{\BC^2} |y_1|^{-2-2\nu_2}|y_2|^{-2-2\nu_1}
(1 + |y_1|^2)^{-2\sigma_1}(1 + |y_2|^2)^{-2\sigma_2}|y_1 - y_2|^{-2-2\nu_3} dy_1dy_2
$$
After one more substitution $y_i \mapsto y_i^{-1},\ i = 1, 2$, 
$$
I'_\BC(\sigma_1,\sigma_2,\sigma_3) = \tI_\BC(\sigma_1,\sigma_2,\sigma_3)
$$
(note that $d(y^{-1}) = dy/|y|^4$). 

Passing to polar coordinates we get 
$$
\int_\BC \frac{dx}{(1 + |x|^2)^2} = \pi,
\eqno{(5.1.6)}
$$
cf. Lemma 2.6. 

Thus
$$
I_\BC(\sigma_1,\sigma_2,\sigma_3) = \pi \tI_\BC(\sigma_1,\sigma_2,\sigma_3)
\eqno{(5.1.6)}
$$
which establishes (1.1.3). $\square$ 

{\bf 5.2.} To establish (1.1.4) it remains to compute the integral (1.1.2). Introduce new parameters 
$a_i = - 1 - \nu_i,\ i = 1, 2, 3.$ In terms of them (1.1.2) becomes
$$
I_{\BC;2}(a_1,a_2,a_3) := 
\int_{\BC^2} (1 + |z_1|^2)^{-2-a_2-a_3}(1 + |z_2|^2)^{-2-a_1-a_3}|z_1 - z_2|^{2a_3} dz_1dz_2.  
$$
In view of 5.1, (1.1.4) is equivalent to

{\bf 5.3. Theorem.} {\it For $\Re a_1, \Re a_2, \Re a_3$ sufficiently large 
$$
I_{\BC;2}(a_1,a_2,a_3) =  
$$
$$
= \pi^2\frac{\Gamma(a_1+a_2+a_3+2)\Gamma(a_1+1)\Gamma(a_2+1)\Gamma(a_3+1)}
{\Gamma(a_1+a_2+2)\Gamma(a_1+a_3+2)\Gamma(a_2+a_3+2)}
\eqno{(5.3.1)}
$$}   

{\bf Proof.} Let us make a change of variables $z_k = r_ke^{i\phi_k}$. If $z_k = dx_k + iy_k$ then 
$dz_k = dx_kdy_k = r_kd\phi_k$. 

So we get:
$$
I_{\BC;3}(a_1,a_2,a_3) = \int_{\BR_+^2} (1 + r_1^2)^{-2-a_2-a_3}(1 + r_2^2)^{-2-a_1-a_3}r_1r_2
$$
$$
\biggl(\int_{[0,2\pi]^2}  |r_1e^{i\phi_1} - r_2e^{i\phi_2}|^{2a_3} d\phi_1d\phi_2\biggr)
dr_1dr_2.
$$
We have
$$ 
|r_1e^{i\phi_1} - r_2e^{i\phi_2}|^2 = (r_1e^{i\phi_1} - r_2e^{i\phi_2})(r_1e^{-i\phi_1} - r_2e^{-i\phi_2}) = 
$$
$$
(r_1 - r_2 e^{i(\phi_2 - \phi_1)})(r_1 - r_2 e^{-i(\phi_2 - \phi_1)})
$$
Suppose that $a_3\in \BN$. 
Note that 
$$
\frac{1}{(2\pi)^2}\int_{[0,2\pi]^2} (r_1 - r_2 e^{i(\phi_2 - \phi_1)})^a(r_1 - r_2 e^{-i(\phi_2 - \phi_1)})^a
d\phi_1d\phi_2 = 
$$
$$
\CT_z (r_1 - r_2 z_1/z_2)^a(r_1 - r_2 z_2/z_1)^a
$$

Let us introduce a polynomial 
$$
\phi_a(r_1,r_2) := \CT_z (r_1 - r_2 z_1/z_2)^a(r_1 - r_2 z_2/z_1)^a = \sum_{i=0}^a \binom{a}{i}^2 r_1^{2i}r_2^{2a - 2i}
$$
Note that 
$$
\phi_a(1, 1) = \binom{2a}{a}
$$
(this is the simplest case of the Dyson identity). 

We get
$$
I_{\BC;2}(a_1,a_2,a_3) = \sum_{i=0}^{a_3} \binom{a_3}{i}^2 
\int_{\BR_+^2} (1 + r_1^2)^{-2-a_2-a_3}(1 + r_2^2)^{-2-a_1-a_3}r_1^{1 + 2i}r_2^{1 + 2a_3 - 2i} dr_1 dr_2
$$
We have 
$$
I(a,b):= \int_0^\infty (1 + r^2)^b r^{1 + 2a} dr = 
$$
($u = r^2$)
$$
\frac{1}{2}\int_0^\infty (1 + u)^b u^{a} du = 
$$
($u = v/(1 - v)$) 
$$
\frac{1}{2}\int_0^1 (1 - v)^{- b - a - 2}v^a dv = \frac{1}{2}B(a+1, - a - b - 1) = 
\frac{1}{2}\frac{\Gamma(a+1)\Gamma(- a - b - 1)}{\Gamma(- b)}
$$
It follows:
$$
\int_{\BR_+^2} (1 + r_1^2)^{-2-a_2-a_3}(1 + r_2^2)^{-2-a_1-a_3}r_1^{1 + 2i}r_2^{1 + 2a_3 - 2i} dr_1 dr_2 = 
$$
$$
\frac{1}{4}\frac{\Gamma(i+1)\Gamma(1 - i + a_2 + a_3)\Gamma(a_3-i+1)\Gamma(1 + i + a_1)}
{\Gamma(a_2 + a_3 + 2)\Gamma(a_1 + a_3 + 2)} 
$$
So
$$
I_{\BC;2}(a_1,a_2,a_3) = \frac{\pi^2\Gamma(a_3+1)^2}{\Gamma(a_2 + a_3 + 2)\Gamma(a_1 + a_3 + 2)}
$$
$$
\sum_{i=0}^{a_3}\frac{\Gamma(1 - i + a_2 + a_3)\Gamma(1 + i + a_1)}
{\Gamma(i+1)\Gamma(a_3 - i + 1)}
$$
Note that if $a_2, a_3\in \BN$, 
$$
\frac{\Gamma(1 - i + a_2 + a_3)\Gamma(1 + i + a_1)}
{\Gamma(i+1)\Gamma(a_3 - i + 1)} = 
\Gamma(a_1 + 1)\Gamma(a_2 + 1)\binom{a_1 + i}{a_1}\binom{a_2 + a_3 - i}{a_2}
$$
Thus
$$
I_{\BC;2}(a_1,a_2,a_3) = \frac{\pi^2\Gamma(a_3+1)^2\Gamma(a_1 + 1)\Gamma(a_2 + 1)}{\Gamma(a_2 + a_3 + 2)\Gamma(a_1 + a_3 + 2)}
$$
$$
\sum_{i=0}^{a_3}\binom{a_1 + i}{a_1}\binom{a_2 + a_3 - i}{a_2}
\eqno{(5.3.2)}
$$
Next we shall use an elementary 

{\bf 5.4. Lemma.} {\it For natural $a_1, a_2, a_3$, 
$$
\sum_{i=0}^{a_3}\binom{a_1 + i}{a_1}\binom{a_2 + a_3 - i}{a_2} = 
\binom{a_1+a_2+a_3+1}{a_3}
$$}

{\bf Proof.} The case $a_2 = 0$ is a consequence of the identity
$$
\binom{a+1}{b+1} = \binom{a}{b+1} + \binom{a}{b}. 
$$
Suppose that 
$a_2\geq 1$ and make the    
induction on $a_3$. The case $a_3 = 0$ is clear. Suppose we have proven the assertion for 
$a_3 \leq n$. We have:
$$
A(n+1):= \sum_{i=0}^{n+1}\binom{a_1 + i}{a_1}\binom{a_2 + n + 1 - i}{a_2} = 
$$
$$
\sum_{i=0}^{n}\binom{a_1 + i}{a_1}\binom{a_2 + n + 1 - i}{a_2} + \binom{a_1 + n + 1}{a_1}
\eqno{(*)}  
$$
But
$$
\binom{a_2 + n + 1 - i}{a_2} = \binom{a_2 + n - i}{a_2} + \binom{a_2 + n - i}{a_2 - 1}
$$
Inserting this into (*) and using the induction hypothesis one finishes the proof. $\square$

\bigskip\bigskip

Combining (5.3.2) with the above lemma we get the assertion of Thm. 5.3 for natural $a_i$. 

On the other hand, one can verify that both sides of (5.3.1) are bounded when for one $i$ $\Re a_i \ra \infty$ and the other 
two $a_j$'s are fixed. So by Carlsson's theorem (cf. Lemma 5.4), the identity (5.3.1) holds true for all $a_i$ with 
$\Re a_i$ sufficiently large (so that the integral converges). 
$\square$

\bigskip\bigskip

\bigskip\bigskip

\centerline{\bf \S 6. Why the real integral is analogous to the complex one}

\bigskip\bigskip

{\bf 6.1.} Consider a change of variables 
$$
x = \tan \alpha,\ y = \tan \beta,\ z = \tan \gamma
$$
We have
$$
\sin(\alpha - \beta) = (\tan\alpha - \tan\beta)\cos\alpha\cos\beta,
$$
$$
1 + x^2 = \cos^{-2}\alpha
$$
and 
$$
d\alpha = \frac{dx}{1 + x^2}
$$

{\bf 6.2.} Consider the integral (1.2.1)
$$
\CI(a,b,c) = \int_{-\pi}^\pi\int_{-\pi}^\pi\int_{-\pi}^\pi 
|\sin(\alpha - \beta)|^{2c}|\sin(\alpha - \gamma)|^{2b}|\sin(\beta - \gamma)|^{2a} d\alpha d\beta d\gamma.  
$$
The function under the integral is $\pi$-periodic with respect to 
each argument. 
It follows that
$$
\CI(a,b,c) = 8\int_{-\pi/2}^{\pi/2}\int_{-\pi/2}^{\pi/2}\int_{-\pi/2}^{\pi/2} 
|\sin(\alpha - \beta)|^{2c}|\sin(\alpha - \gamma)|^{2b}|\sin(\beta - \gamma)|^{2a} d\alpha d\beta d\gamma =  
$$
$$
8\int_{-\infty}^\infty\int_{-\infty}^\infty\int_{-\infty}^\infty 
(1 + x^2)^{-(b+c+1)}(1 + y^2)^{-(a+c+1)}(1 + z^2)^{-(a+b+1)}
$$
$$ 
|x - y|^{2c}|x - z|^{2b}|y - z|^{2a} dx dy dz
$$
We see that this integral is similar to (1.1.1). 

\vspace{2cm}

\newpage

\centerline{\bf References}

\bigskip\bigskip


[B] W.N.Bailey, Generalised hypergeometric series, Cambridge University Press, 1935.  

[BR] J.Bernstein, A.Reznikov, Estimates of automorphic functions, 
{\it Moscow Math. J.} {\bf 4} (2004), 19 - 37. 

[C] L.Carlitz, Note on a $q$-identity, {\it Math. Scand.} {\bf 3} (1955), 281 - 282. 


[Dy] F.J.Dyson, Statistical theory of energy levels of complex systems. I, 
{\it J. Math. Phys.} {\bf 3} (1962), 140 - 156. 

[GGPS] I.Gelfand, M.Graev, I.Pyatetsky-Shapiro, Representation theory and automorphic 
forms.  

[HMW] D.Harlow, J.Maltz, E.Witten, Analytic continuation of Liouville theory, 
arXiv:1108.441.

[Ka] K.W.J.Kadell, A proof of Andrews' $q$-Dyson conjecture for $n=4$, {\it Trans. AMS}, 
{\bf 290} (1985), 127 - 144. 

[K] D.E.Knuth, The art of computer programming, vol. 1, 3rd Edition. 

[S] A.Selberg, Bemerkninger om et multiplet integral, {\it Norsk Matematisk Tidsskrift}, {\bf 26} 
(1944), 71 - 78.

[M] M.L.Mehta, Random matrices, 2nd Edition. 

[Mo] W.G.Morris, II, Constant term identities for finite and affine root systems: 
conjectures and theorems, PhD Thesis, University of Wisconsin - Madison, 1982.  

[T] E.Titchmarsh, The theory of functions. 

[ZZ] A.Zamolodchikov, Al.Zamolodchikov, Conformal bootstrap in Liouville 
field theory, {\it Nucl. Phys.} {\bf B477} (1996), 577 - 605.

[Ze] D.Zeilberger, A proof of the $G_2$ case of Macdonald's root system - Dyson 
conjecture, {\it SIAM J. Math. Anal.} {\bf 18} (1987), 880 - 883. 

\bigskip\bigskip

Institut de Math\'ematiques de Toulouse, Universit\'e Paul Sabatier, 
31062 Toulouse, France

\end{document}